\documentclass[12pt,reqno]{amsart}
\usepackage{amsmath,amssymb,amscd,latexsym,amsthm,mathrsfs,comment}
\usepackage[unicode]{hyperref}
\usepackage{hypbmsec}
\newtheorem{Prop}{Proposition}
\newtheorem{Thm}{Theorem}

\newtheorem{Def}{Definition}
\newtheorem{Cor}{Corollary}
\theoremstyle{remark}

\begin{document}

\title{Counting numbers in multiplicative sets: Landau versus Ramanujan}
\author{Pieter Moree} 
\address{Max-Planck-Institut f\"ur Mathematik, Vivatsgasse 7, D-53111 Bonn, Germany}
\email{moree@mpim-bonn.mpg.de}
\date{\today}
\begin{abstract}
A set $S$ of integers is said to be {\it multiplicative} if for every pair $m$ and $n$ of coprime integers 
we have that $mn$ is in $S$
iff both $m$ and $n$ are in $S$. Both Landau and Ramanujan gave approximations to $S(x)$, the
number of $n\le x$ that are in $S$, for specific choices of $S$. The asymptotical precision of their respective
approaches are being compared and related to Euler-Kronecker constants, a generalization of Euler's
constant $\gamma=0.57721566\ldots$.\\
\indent This paper claims little originality, its aim is to give a survey on the literature related to this theme with an emphasis on the contributions of the author (and his coauthors).
\end{abstract}

\maketitle
\section{Introduction}
\noindent 
To every prime $p$ we associate a set $E(p)$ of positive allowed exponents. Thus $E(p)$ is a subset
of $\mathbb N$. We consider the set $S$ of integers consisting of 1 and
all integers $n$ of the form $n=\prod_i p_i^{e_i}$ with $e_i\in E(p_i)$. Note that this
set is {\it multiplicative}, i.e., if $m$ and $n$ are coprime integers then $mn$ is in $S$
iff both $m$ and $n$ are in $S$. It is easy to see that in this way we obtain all multiplicative
sets of natural numbers. As an example, let us consider the case where $E(p)$ consists of the
positive even integers if $p\equiv 3({\rm mod~}4)$ and $E(p)=\mathbb N$ for the other primes.
The set $S_B$ obtained in this way can be described in another way. By the well-known result
that every positive integer can be written as a sum of two squares iff every prime divisor $p$ of
$n$ of the form $p\equiv 3({\rm mod~}4)$ occurs to an even exponent, we see that $S_B$
is the set of positive integers that can be written as a sum of two integer squares.\\
\indent In this note
we are interested in the counting function associated to $S$, $S(x)$, which counts the
number of $n\le x$ that are in $S$. 
By $\pi_S(x)$ we denote the number of primes $p\le x$
that are in $S$. We will only consider $S$ with the property that $\pi_S(x)$ can be
well-approximated by $\delta \pi(x)$ with $\delta>0$ real and $\pi(x)$ the prime counting
function (thus $\pi(x)=\sum_{p\le x}1$). Recall that the Prime Number Theorem states
that asymptotically $\pi(x)\sim x/\log x$. Gauss as a teenager conjectured that
the logarithmic integral, Li$(x)$, defined as $\int_2^x{dt/\log t}$ gives a much better
approximation to $\pi(x)$. Indeed, it is now known that, for any $r>0$ we have 
$\pi(x)={\rm Li}(x)+O(x\log^{-r}x)$. On the other hand, the result that
$\pi(x)={x/\log x}+O(x\log^{-r}x)$, is false for $r>2$. In this note two types
of approximation of $\pi_S(x)$ by $\delta \pi(x)$ play an important role. We say
$S$ satisfies Condtion A if, asymptotically,
\begin{equation}
\pi_S(x)\sim \delta {x\over \log x}.
\end{equation}
We say that $S$ satisfies Condition B if
there are some fixed positive numbers $\delta$ and $\rho$ such that asymptotically
\begin{equation}
\label{conditionB}
\pi_S(x)=\delta{\rm Li}(x)+O\Big({x\over \log^{2+\rho}x}\Big).
\end{equation}
The following result is a special case of a result
of Wirsing \cite{Wirsing}, with a reformulation following Finch et al.~\cite[p. 2732]{FMS}. As usual
$\Gamma$ will denote the gamma function. By $\chi_S$ we deonte the characteristic function of $S$, that
is we put $\chi_S(n)=1$ if $n$ is in $S$ and zero otherwise.
\begin{Thm}
\label{een}
Let $S$ be a multiplicative set satisfying Condtion A, then
$$S(x)\sim C_0(S) x \log^{\delta-1}x,$$
where 
$$C_0(S)={1\over \Gamma(\delta)}\lim_{P\rightarrow \infty}\prod_{p<P}
\Big(1+{\chi_S(p)\over p}+{\chi_S(p^2)\over p^2}
+\ldots\Big)\Big(1-{1\over p}\Big)^{\delta},$$
converges and hence is positive.
\end{Thm}
In case $S=S_B$ we have 
$\delta=1/2$ by Dirichlet's prime number theoorem for arithmetic progressions. 
Recall that for fixed $r>0$ this theorem states that
$$\pi(x;d,a):=\sum_{p\le x,~p\equiv a({\rm mod~}d)}1={{\rm Li}(x)\over \varphi(d)}+
O\Big({x\over \log^{r}x}\Big).$$
Theorem \ref{een} thus gives that, asymptotically,
$S_B(x)\sim C_0(S_B)x/\sqrt{\log x}$, a result derived in 1908 by Edmund Landau.
Ramanujan, in his first letter to Hardy (1913), wrote in our notation that
\begin{equation}
\label{kleemie}
S_B(x)=C_0(S_B)\int_2^x{dt\over \sqrt{\log t}}+\theta(x),
\end{equation}
with $\theta(x)$ very small. In reply to Hardy's question what `very small' is in this context
Ramanujan wrote back $O(\sqrt{x/\log x})$. 
(For a more detailed account and further references see Moree and Cazaran \cite{MC}.)
Note that by partial integration Ramanujan's claim, if true, implies
the result of Landau. This leads us to the following defintion.
\begin{Def}
Let $S$ be a multiplicative set such that $\pi_S(x)\sim \delta x/\log x$ for some $\delta>0$.
If for all $x$ sufficiently large $$|S(x)- C_0(S) x \log^{\delta-1}x|< 
|S(x)- C_0(S) \int_2^x  \log^{\delta-1}dt|,$$
for every $x$ sufficiently large, we say that the Landau approximation is better than the
Ramanujan approximation. If the reverse inequality holds for every $x$ sufficiently large, we say
that the Ramanujan approximation is better than the Landau approximation.
\end{Def}
We denote the formal Dirichlet series $\sum_{n=1,~n\in S}^{\infty}n^{-s}$ associated
to $S$ by $L_S(s)$. For Re$(s)>1$ it converges. If
\begin{equation}
\label{EK}
\gamma_S:=\lim_{s\rightarrow 1+0}\Big({L_{S}'(s)\over L_{S}(s)}+{\delta\over s-1}\Big)
\end{equation}
exists, we say that $S$ has {\it Euler-Kronecker constant} $\gamma_S$. 
In case $S$ consists of all positive integers we have
$L_S(s)=\zeta(s)$ and it is well known that 
\begin{equation}
\label{gammo}
\lim_{s\rightarrow 1+0}\Big({\zeta'(s)\over \zeta(s)}+{1\over s-1}\Big)=\gamma.
\end{equation}
If the multiplicative set $S$ satisfies
condtion B, then it can be shown that $\gamma_S$ exists. Indeed, we have the following result.
\begin{Thm}
\label{vier1} {\rm \cite{eerstev}.}
If the multiplicative set $S$ satisfies Condition B, then
$$S(x)=C_0(S)x\log^{\delta-1}x\Big(1+(1+o(1)){C_1(S)\over \log x}\Big),\qquad {as}\quad x\to\infty,$$
where $C_1(S)=(1-\delta)(1-\gamma_S)$.
\end{Thm}
\begin{Cor}
Suppose that $S$ is multiplicative and satisfies Condition B. If $\gamma_S<1/2$, then the Ramanujan 
approximation
is asymptotically better than the Landau one. If $\gamma_S>1/2$ it is the other
way around.
\end{Cor}
The corollary follows on noting that by partial integration we have
\begin{equation}
\label{part1}
\int_2^x  \log^{\delta-1}dt=x\log^{\delta-1}x\Big(1+{1-\delta\over \log x}+O\Big({1\over \log^2 x}\Big)\Big).
\end{equation}
On
comparing (\ref{part1}) with Theorem \ref{vier1} we see Ramanujan's claim (\ref{kleemie}), if
true, implies $\gamma_{S_B}=0$.\\
\indent A special, but common case, is where the primes in the set $S$ are, with finitely many exceptions, precisely
those in a finite union of arithmetic progressions, that is, there exists a modulus $d$ and integers
$a_1,\ldots,a_s$ such that for all sufficiently large primes $p$ we have $p\in S$ iff
$p\equiv a_i({\rm mod~}d)$ for some $1\le i\le s$. 
(Indeed, all examples we consider in this paper belong to this special case.)
Under this assumption it can be shown, see Serre \cite{Serre}, that $S(x)$ has an aysmptotic
expansion in the sense of Poincar\'e, that is, for every integer $m\ge 1$ we have
\begin{equation}
\label{starrie}
S(x)=C_0(S)x\log^{\delta-1}x\Big(1+{C_1(S)\over \log x}+{C_2(S)\over \log^2 x}+\ldots+
{C_m(S)\over \log^m x}+O({1\over \log^{m+1}x})\Big),
\end{equation}
where the implicit error term may depend on both $m$ and $S$. In particular 
$S_B(x)$ has an expansion of the form (\ref{starrie}) (see, e.g., 
Hardy \cite[p. 63]{Hardy} for a proof). 

\section{On the numerical evaluation of $\gamma_S$}
We discuss various ways of numerically approximating $\gamma_S$. A few of these approaches involve
a generalization of the von Mangoldt function $\Lambda(n)$ (for more details see Section 2.2 of
\cite{MC}).\\
\indent We define $\Lambda_S(n)$ implicitly by
\begin{equation}
\label{loggie}
-{L_S'(s)\over L_S(s)}=\sum_{n=1}^{\infty}{\Lambda_S(n)\over n^s}.
\end{equation}
As an example let us compute $\Lambda_S(n)$ in case $S=\mathbb N$. Since
$$L_{\mathbb N}(s)=\zeta(s)=\prod_p\Big(1-{1\over p^s}\Big)^{-1},$$
we obtain 
$\log \zeta(s)=-\sum_p \log(1-p^{-s})$ and hence
$$-{L_S'(s)\over L_S(s)}=-{\zeta'(s)\over \zeta(s)}=\sum_p {\log p\over p^s-1}.$$
We infer that $\Lambda_S(n)=\Lambda(n)$, the von Mangoldt function. Recall that
$$\Lambda(n)=
\begin{cases}
\log p & {\rm if~}n=p^e;\\
0 & {\rm otherwise}.
\end{cases}
$$
In case $S$ is a multiplicative semigroup generated by $q_1,q_2,...\ldots$, we have 
$$L_S(s)=\prod_i\Big(1-{1\over {q_i}^s}\Big)^{-1},$$
and we find
$$\Lambda_S(n)=
\begin{cases}
\log q_i & {\rm if~}n=q_i^e;\\
0 & {\rm otherwise}.
\end{cases}
$$
Note that $S_B$ is a multiplicative semigroup. It is generated by $2$, the primes $p\equiv 1({\rm mod~}4)$
and the squares of the primes $p\equiv 3({\rm mod~}4)$.\\
\indent For a more general multiplicative set $\Lambda_S(n)$ can become more difficult in nature as we
will now argue.
We claim that (\ref{loggie}) gives rise to the identity
\begin{equation}
\label{idi1}
\chi_S(n)\log n  =\sum_{d|n}\chi_S({n\over d})\Lambda_S(d).
\end{equation}
In the case $S=\mathbb N$, e.g., we obtain $\log n=\sum_{d|n}\Lambda(d)$.
In order to derive (\ref{idi1}) we use the observation that if $F(s)=\sum f(n)n^{-s}$,
$G(s)=\sum g(n)n^{-s}$ and $F(s)G(s)=H(s)=\sum h(n)n^{-s}$ are formal Dirichlet series, then
$h$ is the Dirichlet convolution of $f$ and $g$, that is $h(n)=(f*g)(n)=\sum_{d|n}f(d)g(n/d)$.
By an argument similar to the one that led us to the von Mangoldt function, one sees that
$\Lambda_S(n)=0$ in case $n$ is not a prime power. Thus we can rewrite (\ref{idi1}) as 
\begin{equation}
\label{idi2}
\chi_S(n)\log n  =\sum_{p^j|n}\chi_S({n\over d})\Lambda_S(d).
\end{equation}
By induction one then finds that $\Lambda_S(p^e)=c_S(p^e)\log p$, where $c_S(p)=\chi_S(p)$ and
$c_S(p^e)$ is defined recursively for $e>1$ by
$$c_S(p^e)=e\chi_S(p^e)-\sum_{j=1}^{e-1}c_s(p^j)\chi_S(p^{e-j}).$$
Also a more closed expression for $\Lambda_S(n)$ can be given (\cite[Proposition 13]{MC}), namely 
$$\Lambda_S(n)=e\log p \sum_{m=1}^e {(-1)^{m-1}\over m} \sum_{k_1+k_2+\ldots+k_m=e}\chi_S(p^{k_1})\chi_S(p^{k_2})
\cdots \chi_S(p^{k_m}),$$
if $n=p^e$ for some $e\ge 1$ and $\Lambda_S(n)=0$ otherwise, or alternatively $\Lambda_S(n)=We\log p$, where
$$W=\sum_{l_1+2l_2+\ldots+el_e=e}{(-1)^{l_1+\ldots+l_e-1}\over l_1+l_2+\ldots+l_e}
\Big({l_1+l_2+\ldots+l_e\over l_1!l_2!\cdots l_e!}\Big)
\chi_S(p)^{l_1}\chi_S(p^2)^{l_2}\cdots \chi_S(p^e)^{l_e},$$
if $n=p^e$ and $\Lambda_S(n)=0$ otherwise, where the $k_i$ run through the natural numbers and the $l_j$ through
the non-negative integers.\\
\indent Now that we can compute $\Lambda_S(n)$ we are ready for some formulae expressing $\gamma_S$ in
terms of this function.

\begin{Thm}
\label{vier}
Suppose that $S$ is a multiplicative set satisfying Condition B. Then
$$\sum_{n\le x}{\Lambda_S(n)\over n}=\delta \log x - \gamma_S+O({1\over \log^{\rho}x}).$$
Moreover, we have
$$\gamma_{S}=-\delta \gamma + \sum_{n=1}^{\infty}{\delta-\Lambda_S(n)\over n}.$$
In case $S$ furthermore is a semigroup generated by $q_1,q_2,\ldots$, then one has 
$$\gamma_S=\lim_{x\rightarrow \infty}\Big(\delta \log x -\sum_{q_i\le x}{\log q_i\over q_i-1}\Big).$$
\end{Thm}
The second formula given in
Theorem \ref{vier} easily follows from the first on invoking the classical definition of $\gamma$:
$$\gamma=\lim_{x\rightarrow \infty}\Big(\sum_{n\le x}{1\over n}-\log x\Big).$$
Theorem \ref{vier} is quite suitable for getting an approximative value of $\gamma_{S}$. The formulae given
there, however, do not allow
one to compute $\gamma_{S}$ with a prescribed numerical precision. For doing that another approach is needed, 
the idea of which is to relate the generating series $L_{S}(s)$ to $\zeta(s)$ and then take
the logarithmic derivative. We illustrate this in Section \ref{SEK}
by showing how $\gamma_{S_D}$ (defined 
in that section) can be computed with high numerical precision.

\section{Non-divisibility of multiplicative arithmetic functions}
Given a multiplicative arithmetic function $f$ taking only integer values, it is an almost immediate obervation
that, with $q$ a prime, the set $S_{f;q}:=\{n:q\nmid f(n)\}$ is multiplicative. 

\subsection{Non-divisibility of Ramanujan's $\tau$}
In his so-called `unpublished' manuscript on the partition and tau functions \cite{BO}, Ramanujan considers
the counting function of $S_{\tau;q}$, where $q\in \{3,5,7,23,691\}$ and $\tau$ is the Ramanujan
$\tau$-function. 
Ramanujan's $\tau$-function is defined as the coefficients of the power series in $q$;
$$\Delta:=q\prod_{m=1}^{\infty}(1-q^m)^{24}=\sum_{n=1}^{\infty}\tau(n)q^n.$$
After setting $q=e^{2\pi i z}$, the function $\Delta(z)$ is the unique normalized cusp form of
weight 12 for the full modular group SL$_2(\mathbb Z)$.
It turns out that $\tau$ is a multiplicative function
and hence the set $S_{\tau;q}$ is multiplicative. Given any such $S_{\tau;q}$, Ramanujan denotes
$\chi_{S_{\tau;q}}(n)$ by $t_n$. He then typically writes: ``It is easy to prove by quite elementary
methods that $\sum_{k=1}^n t_k=o(n)$. It can be shown by transcendental methods that
\begin{equation}
\label{simpelonia}
\sum_{k=1}^n t_k\sim {Cn\over \log^{\delta_q} n};
\end{equation}
and
\begin{equation}
\label{kleemie2}
\sum_{k=1}^n t_k=C\int_2^n{dx\over \log ^{\delta_q} x}+O\Big({n\over \log^r n}\Big),
\end{equation}
where $r$ is any positive number'. Ramanujan claims that $\delta_3=\delta_7=\delta_{23}=1/2$, 
$\delta_5=1/4$ and $\delta_{691}=1/690$. 
Except for $q=5$ and $q=691$ Ramanujan also writes down an Euler
product for $C$. These are correct, except for a minor omission he made in case $q=23$.
\begin{Thm} {\rm (\cite{M}).} For $q\in \{3,5,7,23,691\}$ we have
$\gamma_{S_{\tau;q}}\ne 0$ and thus Ramamnujan's claim {\rm (\ref{kleemie2})} is false for $r>2$.
\end{Thm}
The reader might wonder why this specific small set of $q$. The answer is that in these cases Ramanujan
established easy congruences such as 
$$\tau(n)\equiv \sum_{d|n}d^{11}({\rm mod~}691)$$ that allow one to 
easily describe the non-divisibility of $\tau(n)$ for these $q$. Serre, see \cite{SwD}, has shown
that for every odd prime $q$ a formula of type (\ref{simpelonia}) exists, although no simple congruences
as above exist. This result requires quite sophisticated tools, e.g., the theory of $l$-adic representations.
The question that
arises is whether $\gamma_{S_{\tau;q}}$ exists for every odd $q$ and if yes, to compute it with enough
numerical precision as to determine whether it is zero or not and to be able to tell whether the
Landau or the Ramanujan approximation is better. 

\subsection{Non-divisibility of Euler's totient function $\varphi$}
Spearman and Williams \cite{SW} determined the
asymptotic behaviour of $S_{\varphi;q}(x)$. Here invariants from the cyclotomic field $\mathbb Q(\zeta_q)$ come
into play. The mathematical connection with cyclotomic fields is not very direct in \cite{SW}. However, this
connection can be made and in this way the results of Spearman and Williams can then be rederived in a rather straightforward way, see \cite{FLM, eerstev}. Recall that the Extended Riemann Hypothesis (ERH) says that the
Riemann Hypothesis holds true for every Dirichlet L-series $L(s,\chi)$.
\begin{Thm}{\rm (\cite{FLM}).} \label{eflm}
For $q\le 67$ we have $1/2>\gamma_{S_{\varphi;q}}>0$. For $q>67$ we have 
$\gamma_{S_{\varphi;q}}>1/2$.
Furthermore we have $\gamma_{S_{\varphi;q}}=\gamma+O(\log^2q/\sqrt{q})$, 
unconditionally with an effective constant, $\gamma_{S_{\varphi;q}}=\gamma+O(q^{\epsilon-1})$, unconditionally
with an ineffective constant and  $\gamma_{S_{\varphi;q}}=\gamma+O((\log q)(\log\log q)/q)$ if ERH holds true.
\end{Thm}
The explicit inequalities in this result were 
first proved by the author \cite{eerstev}, who established them assuming ERH. Note that the
result shows that Landau wins over Ramanujan for every prime $q\ge 71$.\\
\indent Given a number field $K$, the Euler-Kronecker constant ${\mathcal EK}_K$ of the number field $K$ is
defined as $${\mathcal EK}_K=\lim_{s \downarrow 1}\Big({\zeta'_K(s)\over \zeta_K(s)}+{1\over s-1}\Big),$$
where $\zeta_K(s)$ denotes the Dedekind zeta-function of $K$. Given a prime $p\ne q$, let $f_p$ the smallest
positive integer such that $p^{f_p}\equiv 1({\rm mod~}q)$. Put
$$S(q)=\sum_{p\ne q,f_p\ge 2}  
{\log p\over p^{f_p}-1}.$$
We have 
\begin{equation}
\label{EK01}
\gamma_{S_{\varphi;q}}= \gamma-{(3-q)\log q\over (q-1)^2(q+1)} -S(q) -
{\mathcal{EK}_{{\mathbb Q}(\zeta_q)}\over q-1}.
\end{equation}
(This is a consequence of Theorem \ref{vier1} and Proposition 2 of Ford et al. \cite{FLM}.)\\
\indent The Euler-Kronecker constants ${\mathcal EK}_K$ and in particular 
$\mathcal{EK}_{{\mathbb Q}(\zeta_q)}$ have been well-studied, see e.g. Ford et al.~\cite{FLM}, Ihara \cite{I} or 
Kumar Murty \cite{KM} for results and references.

\section{Some Euler-Kronecker constants related to binary quadratic forms}
\label{SEK}
Hardy \cite[p. 9, p. 63]{Hardy}  was under the misapprehension that for $S_B$ Landau's approximation is better. However, he based himself
on a computation of his student Geraldine Stanley \cite{Stanley} that turned out to be incorrect.
Shanks proved that
\begin{equation}
\label{geraldine}
\gamma_{S_B}={\gamma\over 2}+{1\over 2}{L'\over L}(1,\chi_{-4})-{\log 2\over 2}
-\sum_{p\equiv 3({\rm mod~}4)}{\log p\over p^2-1}.
\end{equation}
Various mathematicians independently discovered the result that
$${L'\over L}(1,\chi_{-4})=\log\Big(M(1,\sqrt{2})^2e^{\gamma}/2\Big),$$
where $M(1,\sqrt{2})$ denotes the limiting value of Lagrange's AGM algorithm
$a_{n+1}=(a_n+b_n)/2$, $b_{n+1}=\sqrt{a_n b_n}$ with starting values $a_0=1$ and $b_0=\sqrt{2}$.
Gauss showed (in his diary) that
$${1\over M(1,\sqrt{2})}={2\over \pi}\int_0^1 {dx\over \sqrt{1-x^4}}.$$
The total arclength of the lemniscate $r^2=\cos(2\theta)$ is given by $2l$, where
$L=\pi/M(1,\sqrt{2})$ is the so-called lemniscate constant.\\
\indent Shanks used these formulae to show that 
$\gamma_{S_B}=-0.1638973186345\ldots \ne 0$, thus establishing the
falsity of Ramanujan's claim (\ref{kleemie}). 
Since $\gamma_{S_B}<1/2$, it follows by Corollary 1 that actually the Ramanujan approximation is better.

A natural question is to determine the primitive binary quadratic forms $f(X,Y)=aX^2+bXY+cY^2$
of negative discriminant for which the integers represented form a multiplicative set. 
This does not seem to be known. However, in the more restrictive case where we require the multiplicative
set to be also a semigroup the answer is known, see Earnest and Fitzgerald \cite{earnest}.
\begin{Thm}
The value set of a positive definite integral binary quadratic form forms a semigroup if and only if it is in the 
principal class, i.e. represents 1, or has order 3 (under Gauss composition).
\end{Thm}
In the former case, the set of represented integers is just the set of norms from the order 
${\mathfrak O}_D$, which is multiplicative.  In the latter case, the smallest example are the forms of 
discriminant -23, for which the class group is cyclic of order 3: the primes $p$ are partitioned into those of the form $X^2 - XY + 6Y^2$ and those of the form $2X^2 \pm XY + 3Y^2$.

Although the integers represented by $f(X,Y)$ do not in general form a multiplicative set, the associated set
$I_f$ of integers represented by $f$, always satisfies the same type of asymptotic, namely we have
$$I_f(x)\sim C_f{x\over \sqrt{\log x}}.$$  
This result is due to Paul Bernays \cite{Bernays}, of fame in logic, who did his PhD thesis with Landau. Since his work
was not published in a mathematical journal it got forgotten and later partially rediscovered by mathematicians
such as James and Pall. For a brief description of the proof approach of Bernays see Brink et al. \cite{Brink}. 

We like to point out that in general the estimate
$$I_f(x)=C_f{x\over \sqrt{\log x}}\Big(1+(1+o(1)){C'_f\over \log x}\Big)$$  
does not hold. For example, for $f(X,Y)=X^2+14Y^2$, see Shanks and Schmid \cite{SS}.

Bernays did not compute $C_f$, this was only done much later and required the combined effort of various
mathematicians. The author and Osburn \cite{mos} combined these results to show that of all the two dimensional lattices
of covolume 1, the hexagonal lattice has the fewest distances. Earlier Conway and Sloane \cite{CS} had identified the
lattices with fewest distances in dimensions 3 to 8, also relying on the work of many other mathematicians. 

In the special case where $f=X^2+nY^2$, a remark in a paper of Shanks seemed to suggest that he thought $C_f$
would be maximal in case $n=2$. However, the maximum does not occur for $n=2$, see Brink et al. \cite{Brink}.

In estimating $I_f(x)$, the first step is to count $B_D(x)$. 
Given a discriminant $D\le -3$ we let $B_D(x)$ count the number of integers $n\le x$ that are coprime to
$D$ and can be represented by some primitive quadratic integral form of discriminant $D$. The integers so
represented are known, see e.g. James \cite{James}, to form a multiplicative semigroup, $S_D$,  generated by the 
primes $p$ with $({D\over p})=1$ and the squares 
of the primes $q$ with $({D\over q})=-1$. James \cite{James} showed that we have
$$B_D(x)=C(S_D){x\over \sqrt{\log x}}+O({x\over \log x}).$$
An easier proof, following closely the ideas employed by Rieger \cite{Rieger}, was given by Williams \cite{Williams}.
The set of primes in $S_D$ has density $\delta=1/2$. By the law of quadratic reciprocity the set of primes
$p$ satisfying $({D\over p})=1$ is, with finitely many exceptions, precisely a union of arithmetic progressions.
It thus follows that Condition B is satisfied and, moreover, that for every integer $m\ge 1$, 
we have an expansion of the form
$$B_D(x)=C(S_D){x\over \sqrt{\log x}}\big(1+{b_1\over \log x}+{b_2\over \log^2 x}+
\cdots +O({1\over \log^m x})\Big).$$
By Theorem \ref{vier1} and Theorem \ref{vier} we infer that $b_1=(1-\gamma_{S_D})/2$, with
$$\gamma_{S_D}=\lim_{x\rightarrow \infty}\Big({\log x\over 2}-\sum_{p\le x,~({D\over p})=1}{\log p\over p-1}\Big)
-\sum_{({D\over p})=-1}{2\log p\over p^2-1}.$$
As remarked earlier, in order to compute $\gamma_{S_D}$ with some numerical
precision the above formula is not suitable and another route has to be taken.
\begin{Prop} 
\label{expressie} {\rm (\cite{James}.)}
We have, for Re$(s)>1$,
$$L_{S_D}(s)^2=\zeta(s)L(s,\chi_D)\prod_{({D\over p})=-1}(1-p^{-2s})^{-1}\prod_{p|D}(1-p^{-s}).$$
\end{Prop}
\noindent {\it Proof}. On noting that
$$L_{S_D}(s)=\prod_{({D\over p})=1}(1-p^{-s})^{-1}\prod_{({D\over p})=-1}(1-p^{-2s})^{-1},$$
and
$$L(s,\chi_D)=\prod_{({D\over p})=1}(1-p^{-s})^{-1}\prod_{({D\over p})=-1}(1+p^{-s})^{-1},$$
the proof follows on comparing Euler factors on both sides. \qed
\begin{Prop} 
\label{2gamma}
We have
$$2\gamma_{S_D}=\gamma+{L'\over L}(1,\chi_D)-\sum_{({D\over p})=-1}{2\log p\over p^2-1}+
\sum_{p|D}{\log p\over p-1}.$$
\end{Prop}
\noindent {\it Proof}. Follows on logarithmically differentiating the expression for $L_{S_D}(s)^2$ given
in Proposition \ref{expressie}, invoking (\ref{gammo}) and recalling that $L(1,\chi_D)\ne 0$. \qed\\

The latter result together with $b_1=(1-\gamma_{S_D})/2$ leads to a formula first proved by Heupel \cite{Heupel}
in a different way. 

The first sum appearing in Proposition \ref{2gamma} can be evaluated with high numerical precision by using
the identity
\begin{equation}
\label{idie}
\sum_{({D\over p})=-1}{2\log p\over p^2-1}=\sum_{k=1}^{\infty}\Big({L'\over L}(2^k,\chi_{D})-{\zeta'\over \zeta}(2^k)-\sum_{p|D}{\log p\over p^{2^k}-1}\Big).
\end{equation}

This identity in case $D=-3$ was established in \cite[p. 436]{M2}. The proof given there is easily 
generalized. An alternative proof follows on combining Proposition \ref{55} with Proposition \ref{56}.
\begin{Prop}
\label{55}
We have
$$\sum_p {({D\over p})\log p\over p-1}=-{L'\over L}(1,\chi_{D})+
\sum_{k=1}^{\infty}\Big(-{L'\over L}(2^k,\chi_{D})+{\zeta'\over \zeta}(2^k)+\sum_{p|D}{\log p\over p^{2^k}-1}\Big).$$
\end{Prop}
{\it Proof}. This is Lemma 12 in Cilleruelo \cite{C}. \qed
\begin{Prop}
\label{56}
We have
$$-\sum_p {({D\over p})\log p\over p-1}={L'\over L}(1,\chi_{D})+\sum_{({D\over p})=-1}{2\log p\over p^2-1}.$$
\end{Prop}
{\it Proof}. Put $G_d(s)=\prod_p(1-p^{-s})^{(D/p)}$.
We have $${1\over G_d(s)}=L(s,\chi_{D})\prod_{({D\over p})=-1}(1-p^{-2s}).$$ The result then follows
on logarithmic differentiation of both sides of
the identity and the fact that $L(1,\chi_{D})\ne 0$. \qed\\

The terms in (\ref{idie}) can be calculated with MAGMA with high precision and
the series involved converge very fast. Cilleruelo \cite{C} claims
that 
$$\sum_{k=1}^{\infty} {L'\over L}(2^k,\chi_D)=\sum_{k=1}^6 {L'\over L}(2^k,\chi_D)+{\rm Error},~|{\rm Error}|\le 10^{-40}.$$

We will now rederive Shanks' result (\ref{geraldine}). Since there is only one primitive quadratic form
of discriminant -4, we see that $S_{-4}$ is precisely the set of odd integers that can be written as a sum
of two squares. If $m$ is an odd integer that can be written as a sum of two squares, then so can
$2^em$ with $e\ge 0$ arbitrary.
It follows that $L_{S_B}(s)=(1-2^{-s})^{-1}L_{S_{-4}}(s)$ and hence $\gamma_{S_B}=\gamma_{S_{-4}}-\log 2$. On invoking
Proposition \ref{2gamma} one then finds the identity  (\ref{geraldine}).

\section{Integers composed only of primes in a prescribed arithmetic progession}
Consider an arithmetic progression having infinitely many primes in it, that is consider
the progression $a,a+d,a+2d,\ldots$ with $a$ and $d$ coprime.
Let $S'_{d;a}$ be the multiplicative set of integers composed only
of primes $p\equiv a({\rm mod~}d)$. Here we will only consider the simple case where
$a=1$ and $d=q$ is a prime number. This problem is very closely related to that in Section 3.2.
One has 
$L_{S'_{\varphi;q}}(s)=(1+q^{-s})\prod_{p\equiv 1({\rm mod~}q)\atop p\ne q}(1-p^{-s})^{-1}$.
Since $L_{S'_{q;1}}(s)=\prod_{p\equiv 1({\rm mod~}q)}(1-p^{-s})^{-1}$, we then infer that
$$L_{S'_{\varphi;q}}(s)L_{S'_{q;1}}(s)=\zeta(s)(1-q^{-2s})$$
and hence
\begin{eqnarray}
\label{bloep}
\gamma_{S'_{q;1}} & = &\gamma-\gamma_{S'_{\varphi;q}}+{2\log q\over q^2-1}\\
& = & {\log q\over (q-1)^2} +S(q) +{\mathcal{EK}_{{\mathbb Q}(\zeta_q)}\over q-1},\nonumber
\end{eqnarray}
where the latter equality follows by identity (\ref{EK01}).
By Theorem \ref{eflm}, (\ref{bloep}) and the Table in Ford et al. \cite{FLM}, we then arrive after
some easy analysis at the following
result.
\begin{Thm}
\label{eflm2}
For $q\le 7$ we have $\gamma_{S'_{q;1}}>0.5247$. For $q>7$ we have 
$\gamma_{S'_{q;1}}<0.2862$.
Furthermore we have $\gamma_{S'_{q;1}}=O(\log^2q/\sqrt{q})$, 
unconditionally with an effective constant, $\gamma_{S'_{q;1}}=O(q^{\epsilon-1})$, unconditionally
with an ineffective constant and  $\gamma_{S'_{q;1}}=O((\log q)(\log\log q)/q)$ if ERH holds true.
\end{Thm}

\section{Multiplicative set races}
Given two multiplicative sets $S_1$ and $S_2$, one can wonder whether for every $x\ge 0$ we have
$S_1(x)\ge S_2(x)$. We give an example showing that this question is not as far-fetched as one might
think at first sight. Schmutz Schaller \cite[p. 201]{PSS}, motivated by
considerations from hyperbolic geometry, conjectured that the hexagonal lattice is better
than the square lattice, by which he means that $S_B(x)\ge S_H(x)$ for every $x$, where $S_H$ is the set of squared distances
occurring in the hexagonal lattices, that is the integers represented by 
the quadratic form $X^2+XY+Y^2$. It is well-known that
the numbers represented by this form are the integers generated by the primes $p\equiv 1({\rm mod~}3)$, 3 and
the numbers $p^2$ with $p\equiv 2({\rm mod~}3)$. Thus $S_H$ is a multplicative set. If $0<h_1<h_2<...$ are
the elements in ascending order in $S_H$ and $0<q_1<q_2<\ldots $ the elements in ascending order in $S_B$, then
the conjecture can also be formulated (as Schmutz Schaller did) as $q_j\le h_j$ for every $j\ge 1$. Asymptotically
one easily finds that 
$$S_B(x)\sim C_0(S_B){x\over \sqrt{\log x}},~S_H(x)\sim C_0(S_H){x\over \sqrt{\log x}},$$
with $C_0(S_B)\approx 0.764$ the Landau-Ramanujan constant (see Finch \cite[Section 2.3]{Finch})
and $C_0(S_H)\approx 0.639\cdots$. It is thus
clear that asymptotically the conjecture holds true. However, if one wishes to make the above estimates effective, 
matters become much more complicated. Nonetheless, the author, with computational
help of H. te Riele, managed to establish the conjecture of Schmutz Schaller.
\begin{Thm} \label{hoho} {\rm \cite{MteR}.}
If $S_B$ races against $S_H$, $S_B$ is permanently ahead, that is, we have
$S_B(x)\ge S_H(x)$ for every $x\ge 0$.
\end{Thm}
Many of the ideas used to establish the above result were first developed in \cite{M2}. There some other
multiplicative set races where considered. 
Given coprime positive integers $a$ and $d$, let $S'_{d;a}$ be the multiplicative set of integers composed only
of primes $p\equiv a({\rm mod~}d)$. The author established the following result as a precursor to
Theorem \ref{hoho}.
\begin{Thm} {\rm \cite{M2}.}
\label{haha}
For every $x\ge 0$ we have $S'_{3;2}(x)\ge S'_{3;1}(x)$, $S'_{4;3}(x)\ge S'_{3;1}(x)$, $S'_{3;2}(x)\ge S'_{4;1}(x)$ and
$S'_{4;3}\ge S'_{4;1}(x)$.
\end{Thm}
We like to point out that in every race mentioned in the latter result, the associated prime number races have no ultimate winner. For example, already Littlewood \cite{Littlew} in 1914 showed that 
$\pi_{S'_{3;2}}(x)-\pi_{S'_{3;1}}(x)$ has infinitely many sign changes. 
Note that trivially if $\pi_{S'_{d_1;a_1}}(x)\ge \pi_{S'_{d_2;a_2}}(x)$ for every $x\ge 0$, then 
$S'_{d_1;a_1}(x)\ge S'_{d_2;a_2}(x)$ for every $x\ge 0$.

See Granville and Martin \cite{GM} for a nice introduction to prime
number races. 
\section{Exercises}
\noindent {\tt Exercise 1}. The non-hypotenuse numberss $n=1,2,3,4,6,7,8,9,11,12,14,16,\ldots$ are those natural numbers
for which there is no solution of $n^2=u^2+v^2$ with $u>v>0$ integers. The set $S_{NH}$ of non-hypotenuse numbers
forms a multiplicative set that is generated by 2 and all the primes $p\equiv 3({\rm mod~}4)$. 
Show that $L_{NH}(s)=L_{S_B}(s)/L(s,\chi_{-4})$ and hence
$$2\gamma_{NH}=2\gamma_{S_B}-2{L'\over L}(1,\chi_{-4})=\gamma-\log 2 + \sum_{p>2}{({-1\over p})\log p\over p-1}.$$
{\tt Remark}. Put $f(x)=X^2+1$. Cilleruelo \cite{C} showed that, as $n$ tends to infinity,
$$\log {\rm l.c.m.} (f(1),\ldots,f(n))=n\log n +Jn+o(n),$$
with
$$J=\gamma-1-{\log 2\over 2}-\sum_{p>3}{({-1\over p})\log p\over p-1}=-0.0662756342\ldots$$
We have $J=2\gamma-1-{3\over 2}\log 2-2\gamma_{NH}$.\\
\indent Recently the error term $o(n)$ has been improved by 
Ru\'e et al. \cite{madrid} to 
$$O_{\epsilon}\big({n\over \log^{4/9-\epsilon}n}\big),$$ with
$\epsilon>0$.\\
\vfil\eject
\noindent {\tt Exercise 2}. Let $S'_D$ be the semigroup generated by the primes $p$ with $({D\over p})=-1$. It is easy
to see that $L_{S'_D}(s)^2=L_{S_D}(s)^2L(s,\chi_D)^{-2}$ and hence, by Proposition \ref{2gamma}, we obtain
\begin{eqnarray*}
2\gamma_{S'_D}&=& 2\gamma_{S_D}-2{L'\over L}(1,\chi_D)\\
&=& \gamma -{L'\over L}(1,\chi_D)-\sum_{({D\over p})=1}{2\log p\over p^2-1}+\sum_{p|D}{\log p\over p-1}\\
&=&\gamma+\sum_p{({D\over p})\log p\over p-1}+\sum_{p|D}{\log p\over p-1}.
\end{eqnarray*}

\centerline{{\bf Table :} Overview of Euler-Kronecker constants discussed in this paper}
$~~$\\
\begin{center}
\begin{tabular}{|c|c|c|c|}\hline
set & $\gamma_{\rm set} $ & winner & reference \\ \hline \hline
$n=a^2+b^2$ & $-0.1638\ldots$ & Ramanujan & \cite{Sh} \\ \hline
non-hypotenuse & $-0.4095\ldots$ & Ramanujan & \cite{Sh2} \\ \hline
$3\nmid \tau$ & $+0.5349\ldots$ & Landau & \cite{M} \\ \hline
$5\nmid \tau$ & $+0.3995\ldots$ & Ramanujan & \cite{M} \\ \hline
$7\nmid \tau$ & $+0.2316\ldots$ & Ramanujan & \cite{M}  \\ \hline
$23\nmid \tau$ & $+0.2166\ldots$ & Ramanujan & \cite{M}  \\ \hline
$691\nmid \tau$ & $+0.5717\ldots$ & Landau & \cite{M} \\ \hline
$q\nmid \varphi$, $q\le 67$ & $<0.4977$ & Ramanujan & \cite{FLM} \\ \hline
$q\nmid \varphi$, $q\ge 71$ & $>0.5023$ & Landau & \cite{FLM} \\ \hline
$S'_{q;1}$, $q\le 7$ & $>0.5247$ & Landau & Theorem \ref{eflm2} \\ \hline
$S'_{q;1}$, $q>7$ & $<0.2862$ & Ramanujan & Theorem \ref{eflm2} \\ \hline
\end{tabular}
\end{center}
$~~$\\
$~~$\\
\noindent {\tt Acknowledgement}. I like to thank Andrew Earnest and John Voight for helpful 
information concerning 
qudaratic forms having a value set that is multiplicative, and Ana Zumalac\'arregui for sending me \cite{madrid}.

\end{document}